\documentclass{article}
\usepackage{amsmath}
\usepackage{amssymb}
\usepackage{amsfonts}
\usepackage{amscd}
\def\Q{\mathbf Q}
\def\R{\mathbf R}
\def\HH{\mathfrak H}
\def\N{\mathfrak N}
\def\H{\mathcal H}
\def\T{\mathcal T}
\def\ep{\epsilon}
\def\zt{\zeta}
\def\De{\Delta}
\def\Sy{\operatorname{Sym}}
\def\Sym{\operatorname{Symm}}
\def\QS{\operatorname{QSym}}
\def\NS{\operatorname{NSym}}
\def\Qxy{\Q\langle x,y\rangle}
\begin{document}
\title{Quasi-symmetric functions, multiple zeta values, and rooted trees}
\author{Michael E. Hoffman}
\date{June 1, 2006}
\maketitle
\noindent 
My first talk was about the algebra of multiple zeta values,
and the second about Hopf algebras of rooted trees.  A thread that connects
the two is the Hopf algebra $\QS$ of quasi-symmetric functions.  First
defined by Gessel \cite{G}, $\QS$ consists of those formal power series  
$f\in\Q[[t_1,t_2,\dots]]$ (each $t_i$ having degree one), such that
$f$ has bounded degree, and the coefficient in $f$ of
$$
t_{i_1}^{p_1}t_{i_2}^{p_2}\cdots t_{i_k}^{p_k}
$$
equals the coefficient in $f$ of $t_1^{p_1}t_2^{p_2}\cdots t_k^{p_k}$
whenever $i_1<i_2<\dots<i_k$.  As a vector space, $\QS$ is generated by
the monomial quasi-symmetric functions 
$$
M_{p_1 p_2\cdots p_k}=\sum_{i_1<i_2<\dots<i_k}t_{i_1}^{p_1}t_{i_2}^{p_2}\cdots t_{i_k}^{p_k} .
$$
The algebra $\Sy$ of symmetric functions is a proper subalgebra of $\QS$:
for example, $M_{11}$ and $M_{12}+M_{21}$ are symmetric, but $M_{12}$ is not.
\par 
As an algebra, $\QS$ is generated by those monomial symmetric functions 
corresponding to Lyndon words in the positive integers \cite{MR,H1}.
The subalgebra of $\QS^0\subset\QS$ generated by all Lyndon words other than 
$M_1$ has the vector space basis consisting of all monomial symmetric functions
$M_{p_1p_2\cdots p_k}$ with $p_k>1$ (together with $M_{\emptyset}=1$).  
There is a homorphism $\QS^0\to\R$ given by sending each $t_i$ to 
$\frac1{i}$; that is, the monomial quasi-symmetric function 
$M_{p_1\cdots p_k}$ is sent to the multiple zeta value
\begin{equation}
\label{zet}
\zt(p_k,p_{k-1},\dots,p_1)=\sum_{i_1>i_2>\dots>i_k\ge 1}\frac1{i_1^{p_k}i_2^{p_{k-1}}\cdots i_k^{p_1}} .
\end{equation}
In particular, the subalgebra $\Sy^0=\Sy\cap\QS^0$ (which is the subalgebra
of $\Sy$ generated by the power-sum symmetric functions $M_i$ with $i>1$)
has a homomorphic image in the reals generated by the values $\zt(i)$ of the 
Riemann zeta function at integers $i>1$.  It is also convenient to think of
$M_{p_1\cdots p_k}$ as the monomial $x^{p_k-1}y \cdots x^{p_1-1}y$ in the
noncommutative polynomial ring $\Qxy$ (with $\QS^0$ corresponding to
$\HH^0 := \Q1+x\Qxy y$), so that the quantity (\ref{zet}) is the image under a 
homomorphism $\zt:\HH^0\to\R$ of $x^{p_k-1}y\cdots x^{p_1-1}y$.
In fact, it appears that all identities of multiple zeta values follow from 
the interaction between the algebra structure of $\QS$ and a second algebra
structure on $\HH^0$ coming from the shuffle product in $\Qxy$; see, e.g., 
\cite{H2,H4}.
\par
To give $\QS$ the structure of a graded connected Hopf algebra, one defines
a coproduct $\De$ by
$$
\De(M_{p_1\cdots p_k})=\sum_{j=0}^k M_{p_1\cdots p_i}\otimes M_{p_{i+1}\cdots p_k} .
$$
This coproduct makes the power-sum symmetric functions $M_i$ primitive, 
and the elementary symmetric functions $M_{1\cdots 1}$ divided powers.  
Using this Hopf algebra structure,
one can define an action of $\QS$ on $\Qxy$ that makes $\Qxy$ a 
$\QS$-module algebra (see \cite{H2} for details).  In terms of this action 
one can state a result
of Y. Ohno \cite{O} as follows:  for any word $w$ of $\HH^0$ and 
nonnegative integer $i$,
$$
\zt(h_i\cdot w)=\zt(h_i\cdot\tau(w)) .
$$
(Here $\cdot$ denotes the action, 
$h_i$ is the complete symmetric function of degree $i$, 
and $\tau$ is the antiautomorphism of $\Qxy$ that exchanges $x$ and $y$.)
\par
My second talk concerned the relationship between $\QS$ and some
Hopf algebras of trees (or more precisely forests) defined by Kreimer \cite{K} 
and Foissy \cite{F}.
Kreimer's commutative Hopf algebra $\H_K$, which has as its algebra generators 
rooted trees, is the graded dual of the noncommutative Hopf algebra $\T$ of 
rooted trees defined by Grossman and Larson \cite{GL}.  Foissy's 
noncommutative Hopf algebra $\H_F$, which is generated by planar rooted trees, 
is self-dual.
\par
Now $\Sy$ is a self-dual Hopf algebra.  The larger Hopf algebra $\QS$ is
commutative but not cocommutative, and so cannot be self-dual:  its graded
dual is the Hopf algebra $\NS$ of noncommutative symmetric functions
in the sense of Gelfand \it et al. \rm\cite{Get}.  As an algebra, $\NS$ is
the noncommutative polynomial algebra $\Q\langle e_1, e_2,\dots\rangle$,
with $e_i$ in degree $i$, and the $e_i$ are divided powers.  There
is an abelianization homomorphism $\NS\to\Sy$ sending $e_i$ to the 
elementary symmetric function of degree $i$.
\par
The Hopf algebra structure on $\H_K$ is such that the ``ladder'' trees
$\ell_i$ (where $\ell_i$ is the unbranched tree with $i$ vertices) 
are divided powers:  so the map $\phi:\Sy\to\H_K$ sending the $i$th 
elementary symmetric function to $\ell_i$ is a Hopf algebra homomorphism.
In fact, there is a commutative diagram of Hopf algebras
\begin{equation}
\label{d1}
\begin{CD}
\NS @>{\Phi}>> \H_F\\
@VVV @VVV\\
\Sy @>{\phi}>> \H_K
\end{CD}
\end{equation}
where $\Phi$ sends $e_i$ to the unbranched planar rooted tree having
$i$ vertices.  (The map $\H_F\to\H_K$ sends each planar rooted
tree to the corresponding rooted tree, and forgets order in products.)
The commutative diagram (\ref{d1}) dualizes to give
\begin{equation}
\label{d2}
\begin{CD}
\QS @<{\Phi^*}<< \H_F\\
@AAA @AAA\\
\Sy @<{\phi^*}<< \T
\end{CD}
\end{equation}
and the diagram (\ref{d2}) makes it easy to establish some interesting 
properties of the elements of the Hopf algebras involved.
For example, if for a rooted tree $t$ we let $|t|$ be the number of non-root 
vertices of $t$ and $\Sym(t)$ the symmetry group of $t$, then
$$
\kappa_n=\sum_{|t|=n}\frac{t}{|\Sym(t)|}\in\T
$$
can be seen to form a set of divided powers in $\T$, and 
$\phi^*(\kappa_n)=h_n$, the complete symmetric function of degree $n$.
In fact, $\ep_n := (-1)^nS(\kappa_n)$, where $S$ is the antipode in $\T$,
is an element that maps under $\phi^*$ to the $n$th elementary symmetric
function:  further, $n!\ep_n$ is exactly the rooted tree in which $n$
vertices are directly connected to the root.
\par
If we define an operator $\N:\T\to\T$ by $\N(t)=\ell_2\circ t$, where
$\circ$ is the Grossman-Larson product, then we can define coefficients
$n(t;t')$ by
$$
\N^k(t)=\sum_{|t'|=|t|+k}n(t;t')t' .
$$
If $\ell_1=\bullet$ is the tree consisting of just the root vertex,
then $n(\bullet;t)$ is nonzero for every rooted tree $t$:  in the
terminology of Brouder \cite{B}, $n(\bullet;t)$ is the ``tree multiplicity'' 
of $t$.
For a forest $t_1\cdots t_k$ of rooted trees, let $B_+(t_1\cdots t_k)$
be the rooted tree obtained by attaching the root of each $t_i$ to
a new root vertex.  Then using diagram (\ref{d2}) it is easy to see
that
$$
n(\bullet; B_+(\ell_{n_1}\ell_{n_2}\cdots \ell_{n_k}))=
\frac1{m_1!m_2!\cdots}\binom{n_1+\dots+n_k}{n_1\ n_2\ \cdots\ n_k},
$$
where $m_i$ is the number of the $n_j$ equal to $i$.  Cf. equation (1)
of \cite{B}.

\end{document}